\documentclass[11pt]{article}
    \usepackage{amsfonts, amssymb, amsmath, amsthm, latexsym,amscd, xypic}
    \def \c{\mathbb{C}}
    
    \def \r{\mathbb{R}}

    \def\D{\Delta}

    \def \t{\mathfrak{t}}

    \def \Z{\mathcal{Z}}
    \def \A{\mathcal{A}}
    \def \V{\mathcal{V}}
    \def \s{\mathcal{S}}

    \def \.{\cdot}
    \def \dbar{\overline{\partial}}
    \def \deg{\textup{deg}}
    \def \dim{\textup{dim}}

    \def \diag{\textup{diag}}

    \def \Lie{\textup{Lie}}
    \def \ker{\textup{ker}}

    \def \Gr{\textup{Gr}}

    \def\st{such that }
    \def\coh{cohomology}

    \def\gr{\text{Gr}}
    \def\la{\lambda}
    
    \def\({$($}
    \def\){$)$}
   \def\HT{{\mathcal H}_T}
\def\H{{\mathcal H}}
\def\ms{\medskip}
\def\ni{\noindent}
\def\ef{equivariantly formal }

    \def\Black{}
    \theoremstyle{plain}
    \newtheorem{Th}{Theorem}[section]
    
    \newtheorem{Prop}[Th]{Proposition}
    \newtheorem{Cor}[Th]{Corollary}
    \theoremstyle{definition}
    \newtheorem{Ex}[Th]{Example}
    
    \newtheorem{Rem}[Th]{Remark}
\begin{document}
    \title{Vector Fields, Torus Actions  and Equivariant Cohomology}
    \author{Jim Carrell\\ and \\Kiumars Kaveh\\ and \\ Volker Puppe}
\maketitle
    \section{Introduction}
An old result of the first author and David Lieberman
\cite{C-L1,C-L2} says that  if a compact Kaehler manifold $X$
admits a holomorphic vector field $V$  having at least one zero,
then the Dolbeault cohomology algebra $H^*(X, \Omega^*)$ of $X$ is
isomorphic with the bigraded algebra associated to a certain
filtered graded ring. If $V$ has  isolated zeros, this result says
that the \coh \ algebra $H^*(X,\c)$ of $X$ is  isomorphic with the
graded algebra associated to a certain filtration of the
coordinate ring $\c[Z]$ of the scheme $Z$ defined by zero($V$).
Moreover, if $V$ is generated by a $\c^*$-action, then $H^*(X,
\Omega^*)$ is isomorphic with a bigraded ring associated to a
certain filtration of $H^*(Z, \Omega^*)$.

In the original proofs of these results in \cite{C-L1,C-L2}, the
Deligne Degeneracy Criterion is used to show that the spectral
sequence associated to $V$ degenerates at its $E_2$ term.
This shows the existence of a filtration with the desired
properties, but doesn't give any insight into
how to describe it geometrically. There are a number of results
describing the filtration in the isolated case (see
Section 5 and also \cite{Carrell} for some references), but
nothing seems to have been written on the non isolated case,
at least from the standpoint of  the \coh \ algebra.
On the other hand, a result of the third author
on $S^1$-actions in the topological setting (see \cite[p.132]{Puppe2})
implies that the $S^1$-equivariant \coh \  of $X$ plays an
important role if $V$ is generated by a $\c^*$-action
and zero($V$) is finite but nontrivial. The purpose
of this note is to exploit this observation.
In fact, we show that if $V$ is generated by a torus action
which has a nontrivial but otherwise arbitrary fixed point set, then
equivariant Dolbeault  \coh, as recently treated
by S. Lillywhite \cite{Lillywhite} and C. Teleman
\cite{Teleman}, can be used to give direct geometric proofs
of the main results of \cite{C-L1,C-L2}, which avoid using
Deligne degeneracy and which also shed considerable light on
the unexplained filtrations. What is particularly interesting is the picture obtained when
zero$(V)$ has positive dimension.
\Black

    \section{Zeros of Holomorphic Vector Fields and
Cohomology}\label{ZHVF}
    The purpose of this section is to recall the
spectral sequence associated to a  holomorphic vector field
\cite{C-L1}. Let $X$ denote a connected compact Kaehler manifold
    of dimension $n$ with sheaf of holomorphic functions ${\mathcal O}_X$
and sheaves $\Omega_X^p$ of holomorphic $p$-forms for $p>0$.
The contraction operator $i_V:\Omega_X^p\rightarrow
    \Omega_X^{p-1}$ defines the Koszul complex
   \begin{eqnarray*}\label{HVFCS}
    0\rightarrow\Omega_X^n \rightarrow \Omega_X^{n-1}\rightarrow
    \cdots\rightarrow \Omega_X^1 \rightarrow {\mathcal O}_X\rightarrow 0.
    \end{eqnarray*}
In addition, for all $\phi,\omega \in \Omega_X^{*}$,
    \begin{eqnarray*}
    i_V(\phi \wedge\omega )=i_V\phi\wedge\omega +(-1)^p\phi\wedge
    i_V\omega
    \end{eqnarray*}
    if $\phi\in\Omega_X^p$.
Let $A^{p,q}(X)$ denote the smooth forms on $X$ of type $(p,q)$.
The $\dbar$ operator $A^{p,q} \to A^{p,q+1}$ anti-commutes with
$i_V$, so  $(\dbar  - i_V)^2=0.$ Put
\begin{equation}\label{DEFK}
    K^s_X = \bigoplus_{q-p=s} A^{p,q},
    \end{equation}
   and define $D: K^s_X \rightarrow K^{s+1}_X$
   to be  $\dbar  - i_V$. Then $D^2=0$, so we obtain \coh \
    groups $H^s(K_X^*)$. Moreover, $K_X^*$ is a differential graded
    algebra under the exterior product, so the \coh \ groups form a graded
    $\c$-algebra $H^*(K_X^*)$.
    Let $F_\bullet= F_0 \subset F_1 \subset F_2 \subset \cdots \subset
    F_n$ be the filtration of the double complex $A^{*,*}(X)$,
    with $F_i = \bigoplus_{r\leq i}A^{r,*}(X).$
    Since $i_V$ is a derivation, we obtain filtrations $F_\bullet
    H^s(K_X^*)$ for all $s$ \st
    $$F_iH^s(K^*_X)F_j H^t(K^*_X) \subset F_{i+j}H^{s+t}(K^*_X).$$
    Now consider the  spectral sequence
    \begin{equation}\label{SS} E^{-p,q}_1 = H^q(X,\Omega_X^p) \Rightarrow
    H^{q-p}(K^*_X).
    \end{equation}
    \begin{Th}[\cite{C-L1,C-L2}]\label{CL1}
    If $V$ has zeros, then every differential in $(\ref{SS})$
    is zero.
    Consequently $E_1=E_\infty$, and there are ${\c}$-linear
    isomorphisms
    \begin{equation}\label{CAI}
    H^{p+s}\big( X,\Omega_X^p\big)
    \cong F_pH^s (K^*_X)/F_{p-1}H^s(K^*_X),
    \end{equation}
    for every $p\ge 0$ and $s$ which give an
    isomorphism of bigraded $\c$-algebras
    \begin{equation}\label{iso}
    \bigoplus _{p,s} H^{p+s}(X,\Omega_X^p)\cong \bigoplus _{p,s} F_pH^s
(K^*_X)/F_{p-1}H^s(K^*_X).
    \end{equation}
    \end{Th}
    Now let $V$ be generated by a $T=S^1$-action on $X$. Then the
    fixed  point set $Z=X^T$ is smooth and closed
    in $X$ and, as a fixed point scheme, is reduced  as well.
    Moreover, $Z=$zero$(V)$.
    Hence, we may also consider the analogous complex $K^*_Z$ on $Z$ with
differential $D = \dbar$ since $i_V = 0$ on $Z$. Thus
    \begin{equation}\label{quasi}H^s(K^*_Z) = \bigoplus_{q-p=s} H^q
    (Z,\Omega_Z^p).
   \end{equation}
    A basic fact connecting $H^*(X,\Omega^*_X)$ and
    $H^*(Z,\Omega^*_Z)$ is contained in the statement that
    the inclusion $i_Z:Z\to X$ induces a quasi-isomorphism
    of complexes. That is, $i_Z^*$  is an isomorphism of graded
$\c$-algebras
    \begin{equation}\label{THM-isomorphism}
      \bigoplus_s H^s(K_X^*) \cong \bigoplus_s H^s(K_Z^*).
    \end{equation}
    This isomorphism
    is fundamental. There are several proofs of (\ref{THM-isomorphism})
in the literature,
    e.g. \cite{C-S, Fujiki, Feng}, and we will give a
    yet another proof in Section 4. It has several interesting
    consequences, which we now mention. First of all, we can pull back the
    filtration of $H^s(K_X^*)$ to $H^s(K_Z^*)$ to get the isomorphism
    of bigraded $\c$-algebras
    \begin{equation}\label{quasi-iso}
    \bigoplus _{p,s} H^{p+s}(X,\Omega_X^p)\cong \bigoplus _{p,s} F_pH^s
    (K^*_Z)/F_{p-1}H^s(K^*_Z).
    \end{equation}
    However, the geometric description of the filtration on $H^s(K^*_Z)$ is
far
from clear. In fact, finding a direct
    geometric  description of this filtration in terms of
    the decomposition $\bigoplus_{q-p=s} H^q(Z, \Omega^p)$
    is one of the main motivations of this work.
    We can summarize the above discussion in the following corollary
of Theorem \ref{CL1} and (\ref{THM-isomorphism}).
    First, put
    $H^q (X,\Omega_X^p)=H^{p,q}(X)$ and $h^{p,q}(X)=\dim~H^{p,q}(X)$.
    Also, recall the Hodge-Dolbeault
    decomposition $H^m(X,\c)=\bigoplus _{p+q=m} H^{p,q}(X)$ of the \coh \
    of $X$, and note that if $h^{p,q}(X)=0$ for $p\ne q$, then
    $H^*(X,\c)=\bigoplus _p H^{2p}(X,\c)\cong \bigoplus_p H^{p,p}(X)$.
    \begin{Cor} \label{COR1}
    For every $s$,
    \begin{equation}\label{q-p=s}
      \dim~H^s(K_X^*)= \sum_{q-p=s} h^{p,q}(X)=\sum_{q-p=s} h^{p,q}(Z)
    \end{equation} Hence $\dim ~H^s(K_X^*)=0$ if $|s|> \dim Z$ and
    $h^{p,q}(X) = 0$ if $|p-q| > \dim(Z)$.
    Moreover, the  algebra
    $H^*(Z,\Omega_Z^*)$ with the grading of $(\ref{quasi})$
    admits a filtration
    \st $\gr ~ H^*(Z,\Omega_Z^*) \cong H^*(X,\Omega_X^*)$ as bigraded
    $\c$-algebras.
    Thus, if $h^{p,q}(X) = 0$
    whenever $p\ne q$ \(e.g. if $Z$ is finite\), then
    $H^*(Z,\c)=\bigoplus _p H^{2p}(Z,\c)\cong \bigoplus_p H^{p,p}(Z)$,
    and there exists a filtration of $H^*(Z,\c)$ \st
    \begin{equation}\label{pneq} H^*(X,\c)\cong
    \Gr~H^*(Z,\c).\end{equation}
    \end{Cor}
    The last conclusion is mentioned because the condition
    $h^{p,q}(X) = 0$ whenever $p\ne q$ frequently holds for
    spaces with a torus action, e.g. flag varieties,
    toric varieties and, more generally, spherical varieties. See
Section \ref{sec-examples} for more details.
\section{Remarks on Equivariant Cohomology}
    In this section, we will briefly recall the two basic
    definitions of equivariant cohomology due to Borel and Cartan,
    and state a recent result of Lillywhite \cite[\S 5.1]{Lillywhite}
    and Teleman \cite[Theorem 7.3]{Teleman} on
    equivariant Dolbeault cohomology. Suppose $G$ is a compact
    topological group acting on a space $M$. It is well known that there
    exists a contractible space $E_G$ with a free $G$-action.
    The quotient $BG = EG / G$ is called the classifying
    space of $G$. Put
    $$M_G = (M \times EG) / G.$$
    The {\it equivariant cohomology} of
    $M$ over $\c$ is defined to be
    $$H^*_G(M) = H^*(M_G,\c).$$
If $G$ is a compact torus, say $T$, then $H_T^*(pt)=H^*(BT)$
is identified with the polynomial ring $S=\c[\Lie(T)]$.
Thus, $H^*_T(M)$ is an $S$-module (via the natural map $\pi:M_T \to
BT$), and one has the following fundamental fact:
    \begin{Th}[Localization Theorem] \label{thm-localization}
    Suppose the compact torus $T$ acts on a compact manifold $M$. Then
    the kernel as well as the cokernel of the canonical map $$i^*: H^*_T(M)
    \rightarrow H^*_T(M^T),$$ induced by the inclusion $i: M^T
\hookrightarrow M$,
    are torsion modules over $S$.
    Thus if $H^*_T(M)$ is a free module over $S$, then $i^*$ is
    injective. Moreover, $i^*$ becomes an isomorphism after inverting
    elements of a finitely generated multiplicative subset of the
polynomial algebra $S$.
    \end{Th}
    If $H^*_T(M)$ is a free $S$-module, then the action of $T$ on $M$
    is said to be {\it equivariantly formal}.
    Equivalently, $M$ is equivariantly formal if the
    spectral sequence of the fibration $M_T \rightarrow BT$ collapses.
\begin{Rem}\label{ef}
    A Hamiltonian $T$-action on a symplectic manifold
    is \ef. In particular, by a result of Frankel \cite{Frankel}, a
    $\c^*$-action with fixed points on a compact Kaehler manifold
   is equivariantly formal for $T = S^1 \subset \c^*$.
   Moreover, a compact $T$-space whose
   (ordinary) cohomology vanishes in odd degrees is also \ef\ (cf.
   \cite[$\S$14.1]{G-K-M}).
\end{Rem}
   Cartan's construction \cite{Cartan} of equivariant cohomology
   for tori assumes
    $M$ is a smooth manifold, and $T$ acts smoothly on $M$.
Let $\Omega^*(M)$ be the De Rham complex of $\c$-valued forms on $M$.
Define $\Omega^*_T(M)$ to be the complex consisting of of all the
polynomial
maps $f: \Lie(T) \to (\Omega^*(M))^T$. Here the superscript denotes
the $T$-invariants. (Note: this is the same as defining
$\Omega^*_T(M) = (\Omega^*(M) \otimes_\c S)^T$).
    In particular $$\Omega^*_T :=
    \Omega^*_T(pt) = S^T=S.$$
    The grading on $\Omega^*_T(M)$ is defined by $\deg(f)
    = n + 2p$, if $x \mapsto f(x)$ is of degree $p$ in $x$ and $f(x)
    \in \Omega^n(M)$. The differential $$d_T: \Omega^*_T(M)
    \rightarrow \Omega^*_T(M)$$
    of this complex is defined by $$(d_T f)(x) = d(f(x)) - i_{V_x}
    f(x),$$ where $i_{V_x}$ is the contraction with the vector field
    $V_x$ on $M$ generated by $x \in \Lie(T)$. Then $d_T \circ d_T = 0$ and
    $d_T$  increases the degree in $\Omega^*(M)$ by $1$.
    \begin{Th} [Cartan \cite{Cartan}] $H^*_T(M)$ and $H^*(\Omega^*_T(M),
    d_T)$ are isomorphic graded $\c$-algebras.
    \end{Th}
    If $M$ is a complex manifold and $T$ acts via holomorphic
    transformations, a Dolbeault version of
    $T$-equivariant  cohomology is constructed in a similar way.
    For $x \in \Lie(T)$, let $V_x = W_x + \overline{W_x}$ be the splitting
    of the  generating vector field of $x$ into holomorphic and
    anti-holomorphic components. Imitating the  Cartan construction, let
$A_T^{p,*}(M)$ be the complex of all polynomial maps
    $f$ from $\Lie(T)$ to $(A^{p,*}(M))^T$.
    (Note again that this is the same as defining
${A_T}^{p,*}(M) = (A^{p,*}(M) \otimes_\c S)^T$).
    Giving bidegree $(1,1)$ to the generators of $S$ defines a
    bigrading on the algebra $A_T^{*,*}(M) = \bigoplus_{p,q}
    A_T^{p,q}(M)$.
    Define the differential $\dbar_T$ on $A_T^{p,*}(M)$ by
    $$ (\dbar_T f)(x) = \dbar(f(x)) - i_{W_x}
    f(x).$$ This gives a differential on $A^{*,*}_T(M)$ as well.
    The $q$-th cohomology of the complex $(A_T^{p,*}(M),
    \dbar_T)$ is called the $(p,q)$-th {\it equivariant Dolbeault
cohomology}
    of $M$. It is denoted by $H^{p,q}_{T}(M)$. Finally, put
    $$H^m_{T,\dbar}(M)=\bigoplus_{p+q=m}H^{p,q}_{T}(M).$$
    We now state a recent result of Lillywhite and Teleman.
\begin{Th}[Equivariant Hodge Decomposition \cite{Lillywhite, Teleman}]
    \label{thm-equ-hodge}
    Suppose $X$ is a compact Kaehler
    manifold with an \ef $T$-action by
    holomorphic transformations. Then $H^*_{T, \dbar}(X)$
is a free $\c[t]$-module, and there exists an
    isomorphism
    $$H^*_T(X) \cong H^*_{T, \dbar}(X)$$
    of graded $\c$-algebras.
    \end{Th}
Finally, we recall the definition of the equivariant Chern classes
of a vector bundle. Let $E$ be a complex vector bundle over the
Kaehler manifold $X$, equipped with a linear action of $T$ which
lifts the action of $T$ on $X$. The projection map $p: E \to X$
defines a map from $E_T = E \times_G ET$ to $X_T = X \times_T ET$.
This makes $E_T$ a vector bundle over $X_T$. The $r$-th equivariant
Chern class of $E$, denoted by $c^T_r(E)$, is defined to be the
$r$-th Chern class of $E_T$. It is clear that $c^T_r(E) \in
H^{2r}_T(X)$. We will need the following fact: suppose $X$ is connected
and
$T$ acts trivially on it. Let $E$ be a line bundle with a $T$-action as
above.
Let the weight of action of $T$ on each fibre of $E$ be $\omega$. Then
\begin{equation} \label{equ-chern}
c_1^T(E) = -\omega + c_1(E).
\end{equation}
    \section{The Main Result} \label{sec-main}
    As usual, $T$ will be a $1$-dimensional compact torus
    acting on a connected compact Kaehler manifold $X$ of
    dimension $n$ \st $Z=X^T$ is non empty. Let $\V$ be the generating
    vector field of $1 \in \c = \Lie(\c^*)$, and, as before, let $K_X^*$
    denote the total complex of the Koszul complex of the vector field
    $\V$. The purpose of this section is to prove
    the main results in this paper. In particular, we will
    obtain a picture of the filtrations ${F_\bullet} = F_0 \subset F_1
    \subset \cdots \subset F_n$
    of $H^*(K_X^*)$ and $H^*(K_Z^*)$ defined in Section \ref{ZHVF}
    in terms of
    the Dolbeault $T$-equivariant cohomology of $X$.
    First define a map $\Phi_X: H^*_{T,\overline{\partial}}(X)
    \rightarrow H^*(K^*_X)$. Let $f: \t \rightarrow (A^{*,*}(X))^T$ be
an element of $A^{*,*}_T(X)$. It follows
    from the definition (\ref{DEFK}) of $K_X^*$
    that if $f\in A_T^{p,q}(X)$, then $f(1)\in K^{q-p}_X$.
    Therefore, put $$\tilde{\Phi}(f) = f(1).$$
\begin{Prop} \label{commute-differential} $\tilde{\Phi}$ is a
cochain map. That is,
    for $f \in A^{*,*}_T(X)$, we have: $$\tilde{\Phi}(\dbar_T f) =
    D(\tilde{\Phi}(f)).$$
    \end{Prop}
    \begin{proof}
    \begin{eqnarray*}
    \tilde{\Phi}(\dbar_T f) &=& \tilde{\Phi}(\dbar f(x)-i_{V_x}f(x)) \cr
    &=& \dbar f(1) - i_Vf(1) \cr &=& D(f(1)) \cr &=& D(\tilde{\Phi}(f))
    \cr
    \end{eqnarray*}
    Here $V_x$ (respectively $V$) is the generating vector field of
$x \in \t = \Lie(T)$ (respectively $1 \in \t$).
    \end{proof}
Let us  compute  $\dim H^s(K^*_X)$ for all $s$.
Put $\HT^s(X)=\bigoplus _i H_T^{i,i+s}(X)$. Note that from
Theorem \ref{thm-equ-hodge}, $H^*_T(X) = \bigoplus_{s} \H_T^s(X)$.
This gives a new grading on $H^*_T(X)$.
We denote $H_T^*(X)$ together
with this grading by $\H^*_T(X)$. By the above
proposition, $\tilde{\Phi}$ induces a map
   \begin{equation}\label{phisubx} \Phi_{X,s}:\HT^s (X)\to H^s(K_X^*).
\end{equation}
It is not hard to check that this is a $\c$-algebra homomorphism.
For each $p\ge 0$,
we have an edge map $e_{p,s}:F_pH^s(K_X^*) \to H^{p,p+s}(X)$ whose
kernel is $F_{p-1}H^s(K_X^*)$. Let $\pi$ denote the natural map
$\pi: H_T^{p,p+s}(X) \to  H^{p,p+s}(X)$.
Let $f: \c \to \Lie(T), f(t)=\sum_i w_it^i, w_i \in A^{p-i,p+s-i}(X)$
represent
a class in $H^{p,p+s}_T(X)$. It follows from the
definitions that $\Phi_{X,s}(f) = \sum_i w_i$ and $\pi(f) = w_0 =
e_{p,s}(\sum_i w_i)$.
That is, the following diagram is commutative.
$$
\xymatrix{H^{p,p+s}_T(X) \ar[d]^{\pi}
\ar[r]^{\Phi_{X,s}} & F_pH^s(K^*_X) \ar[dl]^{e_{p,s}}\\
H^{p,p+s}(X)}
$$
By equivariant formality, $\pi$ is
surjective. It thus follows that the sequence
$$0\to F_{p-1}H^s(K_X^*) \to F_{p}H^s(K_X^*) \to H^{p,p+s}(X)\to 0$$
is exact. Putting $\H^s(X) = \bigoplus_i H^{i, i+s}(X)$, we conclude
\begin{equation}\label{eq1}
\dim H^s(K_X^*) =\dim \H ^s(X) =\sum _i h^{i,i+s}(X).
\end{equation}
Note that, like before, $H^*(X) = \bigoplus_s \H^s(X)$ gives a
grading on $H^*(X)$. We will denote $H^*(X)$ together with this
grading by $\H^*(X)$. For any $\c$-vector space $V$ and $a\in \c$
, let $V[a]$ denote the $\c[t]$-module structure on $V$ where $t$
acts via multiplication by $a$. Since the action of $T$ on $X$ is
equivariantly formal, by Remark \ref{ef}, we have an exact
sequence of $\c[t]$-modules
\begin{equation} \label{eq-exact-seq}
0\to \c^+[t] \H_T^s(X)\to \HT ^s (X)\to \H ^s (X)[0] \to 0,
\end{equation}
where $\c^+[t]$ denotes the ideal generated by $t$
\cite[$\S$1]{Brion}. Hence
$$\H ^s (X)[0] \cong \HT ^s (X)/ \c^+[t] \H_T^s(X)\cong \H_T^s(X)
\otimes_{\c[t]}\c[0],$$
\Black
and therefore
\begin{equation}\label{eq2}
\dim (\HT ^s (X)\otimes_{\c[t]} \c[0]) =\dim \H ^s (X) = \sum_i
h^{i,i+s}(X).
\end{equation}
We now come to our main theorem.
\begin{Th}\label{main} The following statements hold for any integer $s$.

\ms \ni $(1)$ The $\c[t]$-module map
$$\hat\Phi_{X,s}:\HT^s(X) \otimes_{\c[t]}  \c[1]
    \rightarrow H^s(K^*_X)[1]$$ is a $\c$-linear isomorphism.

\ms \ni $(2)$ The inclusion mapping $i_Z:Z \to X$ induces an
isomorphism $i^*_Z:H^s(K_X^*) \cong H^s(K_Z^*)=\H ^s (Z)$.

\ms \ni $(3)$  In particular, $\sum _i h^{i,i+s}(X)=\sum _i
h^{i,i+s}(Z).$
\end{Th}
\begin{proof}
By the Localization Theorem (Theorem \ref{thm-localization}) and the fact
that $\H_T^s(Z) =H^s(K^*_Z)$,  the map $i_Z^*$ induces an isomorphism
\begin{equation}\label{eq2a}
\HT^s(X)\otimes_{\c[t]}  \c[1] ~{\cong}~ \HT ^s(Z)\otimes_{\c[t]}
\c[1]=H^s(K^*_Z)[1].
\end{equation}
 This implies that $\dim(\HT^s(X)\otimes_{\c[t]}  \c[1])=\dim
H^s(K^*_Z)$. But for any $a\in \c$, $\dim (\HT ^s
(X)\otimes_{\c[t]} \c[a])$ is the rank of $\H^s_T(X)$ as a free
$\c[t]$-module. In particular, $\dim (\HT ^s (X)\otimes_{\c[t]}
\c[0])=\dim (\HT ^s (X)\otimes_{\c[t]} \c[1])$. Hence, (\ref{eq2})
gives ({\it 3}). Since the map in (\ref{eq2a}) is $i_Z^*
\hat\Phi_{X,s}$, it follows that $\hat\Phi_{X,s}$ is injective,
which proves ({\it 1}). The claim in ({\it 2}) is an obvious
consequence.
\end{proof}
\begin{Rem}
Part (1) in Theorem \ref{main} is analogous to the corollary in
\cite[p. 13]{Puppe1}.  The proof of Theorem \ref{main} implies that the
subcomplex of the Koszul complex consisting of $T$-invariant forms
is quasi-isomorphic to the Koszul complex itself. One gets an alternative proof of
 Theorem \ref{main} if one first proves this result directly
(which is similar to the well known result that invariant forms in the deRham
complex determine the deRham cohomology) and then uses the fact that
the equivariant Dolbeault evaluated at $t=1$ is just the invariant Koszul complex.
In this context, the evaluation at $t=1$  is exact, and hence commutes with
homology, whereas the evaluation at $t=0$ is not.
\end{Rem}
Theorem \ref{main} realizes two of the  goals of the paper:
a simple proof that $i^*_Z$ is a quasi-isomorphism, and a
proof of the  isomorphism (\ref{iso}) of Theorem \ref{CL1}
which doesn't  use the Deligne degeneracy criterion.
Let us now comment further on the filtration $F_p H^*(K^*_Z)$.
Let $\hat{\Phi}_X: \H^*_T(X) \otimes_{\c[t]}  \c[1] \to H^*(K_X^*)$ be the
morphism obtained by combining the $\hat{\Phi}_{X,s}$.
Note that $\hat{\Phi}_X$
is a $\c$-algebra isomorphism, but not a isomorphism of graded algebras. However,
$\HT^*(X)\otimes_{\c[t]}  \c[1]$ and $H^*(K^*_X)$ are
both canonically filtered, the former being the  filtration
induced from the grading on $\HT^*(X)$ and the latter being the
filtration introduced earlier.
In fact, if $p\ge 0$, put $F_p\HT^s(X)=
\bigoplus _{i\le p} H^{i,i+s}_T(X)$. If $\Phi_{X,s}$ is the map
defined in (\ref{phisubx}), then, by
definition,
\begin{equation}\label{inclusion}
\Phi_X\big(F_p\HT^s(X) \big) \subset F_pH^s (K_X^*).
\end{equation}
Note that $\Phi_X$ can be described as the map obtained
by composing $\hat{\Phi}_X$ and the
natural map from $\HT^*(X)$ to $\HT^*(X)\otimes_{\c[t]}  \c[1] $
sending $\alpha$ to $\alpha \otimes_{\c[t]}  1$.
We can now give a
geometric description of the filtration of $H^*(K^*_X)$.
Let $\mathcal{R}_X$ denote the
algebra $\H^*_T(X)/\c^+[t] \H_T^*(X)$. Since the ideal $\c^+[t] \H_T^*(X)$ is
homogeneous with respect to the grading of $\H^*_T(X)$, $\mathcal{R}_X$
inherits a grading from
$\H^*_T(X)$.
\begin{Th}\label{filtration}
The mapping
$\Phi_X$ is a surjection of filtered rings.
That is, for all $s$
$$\Phi_{X,s} \big( F_p\HT^s(X) \big) = F_pH^*(K_X^*),$$
and $\mathcal{R}_X$ is isomorphic to both $\H^*(X)$ and $\Gr H^*(K_X^*)$
as graded algebras.
\end{Th}
\begin{proof}
This follows from (\ref{inclusion}) and (1) of Theorem \ref{main}.
\end{proof}
Since the inclusion map $i_Z:Z\to X$ induces a quasi-isomorphism, we
immediately obtain a description of the filtration of $H^*(K^*_Z)$
whose associated graded is $\H^*(X)$.
 \begin{Cor} For each $p\ge 0$,
$$\Phi_Z \circ i^*_Z
\big(\bigoplus_{0\le i \le p} H^{i,i+s}_T(X)\big) = F_pH^s(K_Z^*).$$
\end{Cor}
We will give an example of how to use this result in the next section.
Note also that the natural map
$$\D_p : H_T^{p,p+s}(X) \to F_p H^s (K_X^*) \to
H^p(X,\Omega_X^{p+s})$$ can be described as the $p$-th derivative map
$$\D_p(f)=\dfrac{1}{p!} f^{(p)}(1). $$
\Black
We now use Theorem \ref{main} to prove
a vanishing theorem which extends part of Corollary \ref{COR1}.
    \begin{Th} \label{thm-equ-vanishing}
    If $|p-q| > \dim(Z)$, then $H^{p,q}_{T}(X) = 0$.
    \end{Th}
    \begin{proof} Since $T$ acts trivially on $Z$,
     $$H^*_{T}(Z) = S \otimes_\c H^*(Z).$$
Thus
    $$H^{p,q}_{T}(Z) = \bigoplus_{i \leq \min(p,q)} S^i
    \otimes_\c H^{p-i, q-i}(Z).$$
Since $|p-q|=|(p-i) - (q-i)| > \dim Z$, it follows that
$H^{p,q}_{T}(Z)
    = 0$ as well. From the equivariant Hodge decomposition (Theorem
\ref{thm-equ-hodge})
$H^{p,q}_{T}(X)\subset H_T^{p+q}(X)$. The result now follows from
the localization theorem (Theorem \ref{thm-localization}) since
$i^*\big(H^{p,q}_{T}(X)\big) \subset H^{p,q}_{T}(Z)\subset
H^{p+q}_{T}(Z).$
    \end{proof}
\section{Examples} \label{sec-examples}
\begin{Ex}[Toric varieties] Let $T$ denote the
$n$-dimensional torus $(\c^*)^n$ and let $X$ be a smooth $T$-toric
variety. Let $M$ (respectively $N$) denote the lattice of
characters (respectively $1$-parameter subgroups) of $T$. The real
vector spaces generated by these lattices will be denoted by
$M_\r$ and $N_\r$ respectively. Let $\gamma \in N$ be a
$1$-parameter subgroup in general position, in the sense that the
fixed point set $Z = X^\gamma$ of the $\c^*$-action defined by
$\gamma$ coincides with $X^T$.
Hence $h^{p,q}(X)=0$ if $p\ne q$,
so it follows that $\H^*(X)=H^*(X)$.
Now let $F_0 \subset F_1 \subset \cdots
$ be the filtration in Corollary \ref{COR1} corresponding to the
vector field induced by the $\c^*$-action given by $\gamma$.
Finally, let $\Sigma$ be the fan of $X$ in $N_\r$. Each ample
$T$-equivariant line bundle $L$ on $X$ gives rise to an integral
polytope $\Delta(X, L) \subset M_\r$. This polytope is the image
of $X$ under the moment map of $(X, L)$. It can be seen that
$\Delta(X, L)$ is normal to the fan $\Sigma$, i.e. every facet of
$\Delta$ is orthogonal to a $1$-dimensional cone of $\Sigma$. Let
$\s$ be the semi-group of all rational convex polytopes normal to
$\Sigma$, where the semi-group operation is the Minkowski sum of
polytopes. Let $\V$ denote the complex vector space generated by
the semi-group $\s$. For $\Delta_1$ and $\Delta_2 \in \s$, we
write $\Delta_1 \sim \Delta_2$ if $\Delta_1$ can be identified
with $\Delta_2$ by a translation. This relation $\sim$ extends to
the whole $\V$. The assignment $c_1(L) \mapsto \Delta(X, L)$
extends to give an isomorphism $H^2(X, \c) \cong \V / \sim$. Since
$Z$ is isolated, we have $H^*(Z, \c) = H^0(Z, \c) \cong \c^Z$.
Thus, every element of $H^0(Z, \c)$ can be regarded as a function
from $Z$ to $\c$. Let $\Delta \in \s$. To each fixed point $z \in
Z = X^T$, there correspond a unique vertex $v_z$ of $\Delta$.
Consider the function $f_\Delta: Z \to \c$ defined by
$$f_\Delta(z) = \langle \gamma, v_z \rangle.$$
The following result is proved in \cite{Kaveh}.
\begin{Prop} \label{prop-toric-filtration}
Under the isomorphism in Corollary \ref{COR1}, the function
$f_\Delta$ corresponds to $\Delta$ regarded as an element of $\V /
\sim \,\cong H^2(X, \c)$.
\end{Prop}
For each fixed point $z \in Z$, define a linear function $\tilde{z}$ on
$\V$
by $$\tilde{z}(\Delta) = \langle \gamma, v_z \rangle.$$ Let $\Z =
\{ \tilde{z} \mid z \in Z \}$ and let $\c[\Z]$ denote the coordinate
ring of $\Z$ as a subset of the vector space $\V^*$. It carries a natural
filtration induced by the degree of polynomials. Since $H^*(X, \c)$ is
generated
by $H^2(X, \c)$, from Corollary \ref{COR1} and
Proposition \ref{prop-toric-filtration}
one can see that $$H^*(X, \c) \cong \Gr \c[\Z].$$
The equivariant cohomology $H^*_T(X, \c)$ can be identified with the
algebra $\A$ of all continuous functions on $N_\r$ whose restriction
to each cone of $\Sigma$ is given by a polynomial (conewise
polynomial). Under this identification $H^{2i}_T(X,
\c)$ corresponds to the subspace $\A_i$ of all continuous functions
whose restriction to each cone of maximal dimension is given by a
homogeneous degree $i$ polynomial. Each fixed point $z$ corresponds
to a cone of maximal dimesnion $\sigma_z$ in $\Sigma$.
Let $g \in \A$. Define the function $\tilde{g}: Z \to \c$ by
$$\tilde{g}(z) = g_{\mid \sigma_z}(\gamma).$$
It is shown in \cite{Kaveh} and also follows from
Theorem \ref{filtration} that $g \in \A_i$ if and only if $\tilde{g}
\in F_i$.
\end{Ex}

\begin{Ex}[The flag variety $G/B$]
Let $G$ be a connected reductive group over $\c$, $B$ a Borel
subgroup
and $X=G/B$ the flag variety of $G$. Let
$T$ be a maximal torus  in $B$ and $\t$ and
$\t^*$ respectively be the Lie algebra of $T$ and its dual.
Also let $S=\c[\t]$ be the algebra of polynomials on $\t$.
The torus $T$ acts on
$X$ by multiplication from left, and any
$1$-parameter subgroup of $T$ induces a
$\c^*$-action on $X$. Let $s\in \t$ induce a
$1$-parameter subgroup of $T$, and suppose
$F_\bullet = F_0 \subset F_1 \subset
\cdots$ is the filtration in Corollary \ref{COR1} corresponding to
this $\c^*$-action. Since $h^{p,q}(X) = 0$ for $p \neq q$,
it follows that $\H^*(X)=H^*(X)$. We now
consider the two obvious cases.

\ms \noindent{\bf (a)} The element $s$ is regular. That is, the
isotropy group $W_s$  is trivial. Then $Z = (G/B)^s = (G/B)^T$ and
$Z$ is identified with the Weyl group $W$ of $(G,T)$ by
identifying $w$ and $wB$. Thus $H^*(Z, \c) = \c^Z = \c^W$. We will
first describe the map $H^*_T(G/B) \to H^*(Z, \c)$ obtained by
localizing and setting $t=1$. Now, $H^*_T(X)$ is isomorphic as a
$\c$-algebra to $S \otimes_{S^W} S$ where $S^W$ denotes the
subalgebra of $W$-invariants (see \cite[$\S 2$ Examples]{Brion}).
Since $H^*(G/B, \c)$ is generated by the Chern classes of line
bundles, and such line bundles are always $T$-equivariant, we need
only consider the image of an equivariant Chern class
$c^T_1(L_\lambda)$, where $L_\lambda$ denotes the line bundle
corresponding to a weight $\lambda\in \t^*.$ But it can be shown
that $c^T_1(L_\lambda)=-\sum_{w \in W} 1 \otimes (w\cdot
\lambda)$, and  so $c^T_1(L_\lambda)$ is sent to the element
$f_\la \in H^*(Z, \c)$ defined by the condition
\begin{equation} \label{chern-flag-a}
f_\la(w)=-\langle w\cdot \lambda, s \rangle.
\end{equation}
This coincides with the representative of $c_1(L_\lambda)$ on
$H^0(Z, \c)$ calculated for example in \cite{contemp}. The
filtration $F_i$ of $H^*(Z, \c)$ was originally obtained as
follows. Consider the $\c$-algebra morphism $S \to \c^W$
determined by
  $\la\in \t^* \to f_\la$. Since any function vanishing
on $W\. s$ is mapped to 0, we get a $\c$-algebra  morphism
$\c[W\cdot s]\to \c^W$,
where $\c[W\. s]$ is the ring of polynomials on the finite set $W\.s$.
This mapping turns out to be an isomorphism sending  the canonical
filtration of
$\c[W\cdot s]$ (defined by the fact that
$\c[W\.s]$ is a quotient of $S$) onto the filtration of
$\c^W=H^*(Z, \c)$. Furthermore, since $s$ is regular, it is not hard
to see that  $\gr~\c[W\cdot s]\cong S/I_W^+$,
where $I_W^+$ denotes the ideal in $S$ generated by the
homogeneous elements of $S^W$. Thus, we have the well
known  Borel picture of $H^*(G/B, \c)$.

\ms \noindent {\bf (b)} The element $s$ is non-regular. Let $\Phi$
be the root system of $(G, T)$ and $\Phi_s = \{ \alpha \in \Phi
\mid \alpha(s) = 0 \}$. Put $$\mathfrak{t}_0 = \bigcap_{\alpha \in
\Phi_s} \ker \alpha,$$ and let $T_0 \subset T$ be the
corresponding torus. Finally, let $L$ denote the Levi subgroup $L
= Z_G(T_0)$. For example, if $G = GL(n, \c)$, and $T$ is the
diagonal torus, put $s = \diag(a_1I_{n_1}, a_2I_{n_2}, \ldots,
a_{n_r}I_{n_r})$ where $I_l$ is the $l \times l$ identity matrix,
$n_1 + \cdots + n_r = n$ and $a_i \neq a_j$ when $i \neq j$. Then
$L = GL(n_1, \c) \times \cdots \times GL(n_r, \c)$. Then the Weyl
group $W_L$ of $L$ is the isotropy subgroup of $s$ in $W$. Now
$Z=(G/B)^s = (G/B)^{T_0}$ is a union of the flag varieties of $L$.
More precisely, for $w \in W$, let $Z_w = $ connected component of
$Z$ containing $wB \in (G/B)^T$. One sees that each $Z_w$ is
isomorphic to $L / L\cap B$ and $Z_w = Z_{w'}$, for $w, w'$ in the
same right coset of $W_L$. Thus $$Z  = \bigcup_{w \in
W_L\backslash W} Z_w.$$ Hence $H^*(Z, \c) = \bigoplus_{w \in W_L
\backslash W} H^*(L/L\cap B)$.
To obtain the filtration of $H^*(Z, \c)$, take an element $t \in
\t$ which determines a regular
$1$-parameter subgroup of $T$. Let $\Z \subset \t \oplus \t$ be the
$W$-orbit of $(s,t)$, where $W$ acts on $\t \oplus \t$ diagonally.
One can write
$$\Z = \bigcup_{w \in W_L \backslash W} \Z_w,$$ where $\Z_w = \{
(w^{-1}\. s, w^{-1}u^{-1}\.t) \mid u \in W_L \}$. The elements of $\Z$ are
in
one-to-one correspondence with $(G/B)^T$, and each $\Z_w$
corresponds to the $T$-fixed points in $Z_w$. Let
$\c[\Z]$ (respectively $\c[\Z_w]$) denote the coordinate ring of
$\Z$ (respectively $\Z_w$). From part (a), $H^*(Z_w)
\cong \Gr ~\c[\Z_w]$ for any $w \in W_L
\backslash W$, , where the filtration on $\c[\Z_w]$
is induced by the degree. Hence
$$H^*(Z, \c) \cong \bigoplus_{w \in W_L \backslash W} \Gr ~\c[\Z_w].$$
Put $\A=\bigoplus_{w \in W_L \backslash W} \Gr ~\c[\Z_w].$ The following
shows that the filtration on $\A$ is induced by the natural filtration on
$\c[\Z]$ given by the degree.
\begin{Prop}
An element $(f_w) \in \A$ lies in $F_i$ if and only if there
exists an element $f \in \c[\Z]$ such that $f$ has degree $\leq i$
and restriction of $f$ to $\Z_w$ is a representative for $f_w$ in
$\Gr ~\c[\Z_w]$.
\end{Prop}
\begin{proof}
Since $H^*(G/B, \c)$ is generated by $H^2(G/B, \c)$, the
filtration is generated by $F_1$, i.e. $F_i$ consists of all
polynomials of degree $\leq i$ in the elements of $F_1$. Hence it
is enough to verify the claim for $F_1$. Consider the line bundle
$L_\lambda$ on $G/B$ corresponding to a dominant weight $\lambda$,
and let $L_{\lambda, w}$ be the restriction of $L_\lambda$ to the
small flag variety $Z_w$. Then, for each $w \in W_L \backslash W$,
the weight of the action of $s$ on $L_{\lambda, w}$ is $\langle
\lambda, w^{-1}\.s \rangle$. From (\ref{equ-chern}),
\begin{equation} \label{chern-flag-b} c^s_1(L_{\lambda, w}) = -\langle
\lambda, w^{-1}\.s\rangle + c_1(L_{\lambda, w}),\end{equation} where
$c^s_1$ denotes the equivariant Chern class for the $\c^*$-action
induced by $s$. Then, from Theorem \ref{filtration},
(\ref{chern-flag-a}) and (\ref{chern-flag-b}) it follows that
$c_1(L_{\lambda, w})$ corresponds to the element $(f_{\lambda, w})$
represented by the function
\begin{equation} \label{filtration-flag} (w^{-1}\. s, w^{-1}u^{-1}\.t)
\mapsto -\langle \lambda, w^{-1}\.s\rangle - \langle \lambda,
w^{-1}u^{-1}\.t\rangle .\end{equation} Now let $f_\lambda$ be the
linear function on $\t \oplus \t$ given by
$$f(x, y) = -\lambda(x) - \lambda(y).$$ From (\ref{filtration-flag}),
the restriction of $f_\lambda$ to $\Z_w$ gives a representative for
$f_w \in \Gr~ \c[\Z_w]$. Since $H^2(G/B, \c)$ is spanned by the
$c_1(L_\lambda)$, this proves the Proposition.
\end{proof}
\end{Ex}

    \noindent James B. Carrell, University of British Columbia,
    Vancouver, B.C., Canada \\{\it Email address:} {\sf
carrell@math.ubc.ca}\\
    \noindent Kiumars Kaveh, University of British Columbia, Vancouver,
    B.C., Canada \\{\it Email address:} {\sf kaveh@math.ubc.ca}\\
    \noindent Volker Puppe,
Universit\"{a}t Konstanz, Fakult\"{a}t f\"{u}r Mathematik, Fach D202,
D-78457 Konstanz,
Germany
\\{\it Email address:} {\sf Volker.Puppe@uni-konstanz.de}
\end{document}